\documentclass[11pt,a4paper,reqno]{article}
\usepackage{amsmath}
\usepackage{amsthm}
\usepackage{enumerate}
\usepackage{amssymb}

\usepackage{verbatim}
\usepackage{url}
\usepackage[latin1]{inputenc}

\usepackage{caption}
\usepackage{graphicx,color}
\def\Z{{\mathbb Z}}

\def\be{\begin{equation}}
\def\ee{\end{equation}}
\def\bea{\begin{equation*}}
\def\eea{\end{equation*}}
\def\bal{\begin{aligned}}
\def\eal{\end{aligned}}

\def\Pr{{\mathbb P}}
\DeclareMathOperator{\E}{{\mathbb E}}

\newtheorem{thm}{Theorem}

\theoremstyle{remark}

\newtheorem{preex}[thm]{Example}

\theoremstyle{definition}

\begin{document}

\title{Tertiles and the time constant}
\date{}
\author{Daniel Ahlberg\thanks{Department of Mathematics, Stockholm University. This work was in part supported by the Swedish Research Council (VR) through grant 2016-04442.}}
\maketitle

\begin{abstract}
We consider planar first-passage percolation and show that the time constant can be bounded by multiples of the first and second tertiles of the weight distribution. As a consequence we obtain a counter-example to a problem proposed by Alm and Deijfen.

\end{abstract}

\section{A theorem and a counter-example}

\vspace{-1pt}

In first-passage percolation on the square lattice, weights $\omega_e$ are assigned independently to the edges according to some distribution $F$ on $[0,\infty)$. The resulting weighted graph induces a random (pseudo-)metric $T$ on $\Z^2$ as follows: For all $x,y\in\Z^2$, let
$$
\textstyle{T(x,y):=\inf\big\{\sum_{e\in\pi}\omega_e:\pi\text{ is a self-avoiding path connecting $x$ to $y$}\big\}}.
$$
Let $Y$ denote the minimum of four independent variables distributed as $F$. When $\E[Y]<\infty$ the limit $\mu:=\lim_{n\to\infty}\frac1nT\big((0,0),(n,0)\big)$ exists a.s.\ as a consequence of Kingman's subadditive ergodic theorem. However, as is well-known, the limit exists in probability for all weight distributions. See e.g.~\cite{aufdamhan17} for background and details.

Upper bounds for $\mu$ can be expressed in terms of moments involving $F$; see e.g.~\cite{smywie78}. The moral of this note is that moments are in general poor estimates of $\mu$.
We provide bounds in terms of the first and second tertiles.
A related lower bound has previously been obtained by Cox~\cite{cox80}. Our proof is much shorter.
Let $t_q:=\inf\{t\ge0:F(t)\ge q\}$.

\begin{thm}
For any $F$ the time constant $\mu$ satisfies
$
\tfrac{1}{100}t_{1/3}\,\le\,\mu\,\le\,2t_{2/3}.
$
\end{thm}

In $d\ge2$ dimensions the arguments can be adapted to give $\frac14 t_{1/2d}\le\mu\le dt_{p_c(d)}$, where $p_c(d)\sim1/d$ is the critical probability for oriented percolation on $\Z^d$.

\begin{proof}
The connoisseur will note that the upper bound is immediate from the `flat edge' of Durrett and Liggett~\cite{durlig81}. Spelling things out, let $A_n$ denote the event that there exists a path of length $n$ connecting the origin to a point in $\{(x,y)\in\Z^2:x+y=n,x\ge n/2,y\ge0\}$ 
having total weight at most $nt_{2/3}+M$. Similarly, let $A_n'$ denote the event that there is a path of length $n$ connecting $(n,0)$ to $\{(x,x)\in\Z^2:0\le x\le n/2\}$ 
having total weight at most $nt_{2/3}+M$. Since $2/3$ exceeds the critical probability for oriented percolation on $\Z^2$ (see~\cite{liggett95}), standard results in percolation theory (see~\cite{durlig81}) show that for $M$ large we have $\Pr(A_n)=\Pr(A_n')\ge3/4$, uniformly in $n$.
On the intersection $A_n\cap A_n'$, which occurs with probability at least $1/2$, we have
$T\big((0,0),(n,0)\big)$ bounded by $2nt_{2/3}+2M$,
thus implying that $\mu\le 2t_{2/3}$.


The lower bound is inspired by an argument explored by Smythe and Wierman~\cite{smywie78}, who in turn cite Hammersley. Given $\delta>0$, let $N_n$ denote the number of self-avoiding walks of length $n$ starting at the origin that have fewer than $\delta n$ edges with $\omega_e\ge t_{1/3}$. The number of self-avoiding walks of length $n$ is at most $2.7^n$ for large $n$ (see~\cite[p.~24]{smywie78}). For a given path of length $n$, the number of edges with $\omega_e\ge t_{1/3}$ is binomially distributed with parameters $n$ and $p\ge 2/3$. Let $X$ be binomial with parameters $n$ and $2/3$. For $\beta>0$ Markov's inequality gives
$$
\Pr(X<\delta n)\,=\,\Pr\big(n-X>(1-\delta)n\big)\,\le\, e^{-\beta(1-\delta)n}\E\big[e^{\beta(n-X)}\big]\,=\,\big[\tfrac13e^{\beta\delta}+\tfrac23e^{-\beta(1-\delta)}\big]^n.
$$
Set $\beta=5$ and $\delta=1/100$. Monotonicity of the binomial distribution then gives
$$
\Pr(N_n\ge1)\,\le\,\E[N_n]\,\le\,2.7^n\,\Pr(X<\delta n)\,\le\, 2.7^n\cdot .36^n\,=\,.972^n.
$$
On $\{N_n=0\}$ we have $T\big((0,0),(n,0)\big)\ge \frac{n}{100} t_{1/3}$, so $\mu\ge \frac{1}{100}t_{1/3}$, as required.
\end{proof}

We have noted that for many common distributions $2t_{2/3}$ exceeds the mean of $F$.
Nevertheless, it is curious that one may disregard as much as a third of the mass of a distribution and yet produce general upper and lower bounds on $\mu$.
A more careful analysis should be able to increase the fraction $1/3$, and decrease $2/3$, arbitrarily close to $1/2$.
A conversation with Michael Damron led to the following examples, suggesting that the median cannot be used to obtain general upper nor lower bounds on $\mu$: Let $(F_n)_{n\ge1}$ put mass $1/2$ at $1$, be fully supported on $[0,1]$, and converging weakly to the balanced Bernoulli distribution. Let $(F_n')_{n\ge1}$ put mass $1/2$ at $1$, be fully supported on $[1,\infty)$, with mass $1/2$ diverging in the limit. In all cases the median is $1$. By continuity the time constant for $F_n$ tends to zero as $n\to\infty$. We believe further that the time constant for $F_n'$ tends to infinity with $n$, although this requires an argument.

Based on simulations, it was suggested in~\cite[page~668]{almdei15}, and restated in~\cite[Question~12]{aufdamhan17}, that $\mu\ge\E[Y]$ should hold in great generality, at least when $F$ puts no mass at zero. The above upper bound on $\mu$ provides a fairly general counter-example:
Let $F$ be any distribution on $[0,\infty)$ that puts mass at least $2/3$ on $[0,1]$ and mass at least $1/6$ on $[3888,\infty)$. Then $\mu\le2$, whereas an easy calculation gives $\E[Y]\ge3$.


\end{document}